\documentclass[10pt,amsfonts]{article}
\usepackage{graphicx}
\usepackage{latexsym}
\usepackage{amssymb}
\usepackage{amsmath}
\usepackage{enumerate}
\usepackage{stmaryrd}
\usepackage{layout}
\usepackage{eufrak}
\newtheorem{prop}{Proposition}
\newtheorem{lemma}{Lemma}
\newtheorem{definition}{Definition}

\newtheorem{theorem}{Theorem}
\newtheorem{remark}{Remark}
\newtheorem{example}{Example}

\def\real{{\mathord{{\rm I\kern-2.8pt R}}}}        
\def\inte{{\mathord{{\rm I\kern-2.8pt N}}}}

\def\sZZ{{\rm Z\kern-2.8ptem{}Z}}

\def\z{{\mathchoice
  {\sZZ}
  {\sZZ}
  {\rm Z\kern-0.30em{}Z}
  {\rm Z\kern-0.25em{}Z} }}
\def\sQQ{{\kern 0.27em \vrule height1.45ex width0.03em depth0em
          \kern-0.30em \rm Q}}
\def\qu{{\mathchoice
    {\sQQ}
    {\sQQ}
  {\kern 0.225em \vrule height1.05ex width0.025em depth0em \kern-0.25em \rm Q}
  {\kern 0.180em \vrule height0.78ex width0.020em depth0em \kern-0.20em \rm Q}
        }}
\def\sCC{{\kern 0.27em \vrule height1.45ex width0.03em depth0em
          \kern-0.30em \rm C}}
\def\complex{{\mathchoice
    {\sCC}
    {\sCC}
  {\kern 0.225em \vrule height1.05ex width0.025em depth0em \kern-0.25em \rm C}
  {\kern 0.180em \vrule height0.78ex width0.020em depth0em \kern-0.20em \rm C}
        }}

\newcommand{\ba}{\begin{array}}
\newcommand{\ea}{\end{array}}
\newcommand{\be}{\begin{equation}}
\newcommand{\ee}{\end{equation}}
\newcommand{\bea}{\begin{eqnarray}}
\newcommand{\eea}{\end{eqnarray}}
\newcommand{\beaa}{\begin{eqnarray*}}
\newcommand{\eeaa}{\end{eqnarray*}}

\def\z{\zeta}

\font\tenmath=msbm10 \font\sevenmath=msbm7 \font\fivemath=msbm5
\newfam\mathfam \textfont\mathfam=\tenmath
\scriptfont\mathfam=\sevenmath \scriptscriptfont\mathfam=\fivemath

\def \={{\buildrel {\rm (law)} \over =}}

\def\HH{\EuFrak H}
\def\qed{ \hfill \vrule width.25cm height.25cm depth0cm\smallskip}

\newcommand{\basa}{\begin{assumption}}
\newcommand{\easa}{\end{assumption}}

\newcommand{\bas}{\begin{assum}}
\newcommand{\eas}{\end{assum}}

\newcommand{\ignore}[1]{}
\begin{document}

\renewcommand{\thefootnote}{\fnsymbol{footnote}}

\renewcommand{\thefootnote}{\fnsymbol{footnote}}
\begin{center}
{\large{\bf Cram\'er theorem for Gamma random variables}}\\~\\
Solesne Bourguin\footnote{
SAMM, Universit\'e de Paris 1 Panth\'eon-Sorbonne, 90, rue de Tolbiac, 75634, Paris, France.
Email: {\tt solesne.bourguin@univ-paris1.fr}
}
and
Ciprian A. Tudor\footnote{Laboratoire Paul Painlev\'e, Universit\'e de Lille 1, F-59655 Villeneuve d'Ascq, France. 
Email: {\tt tudor@math.univ-lille1.fr}}\footnote{Associate member of the team Samm, Universit\'e de Paris 1 Panth\'eon-Sorbonne}\\
{\it Universit\'e Paris 1} and {\it Universit\'e Lille 1}\\~\\
\end{center}
{\small \noindent {\bf Abstract:} 
In this paper we discuss the following problem: given a random variable $Z=X+Y$ with Gamma law such that $X$ and $Y$ are independent, we want to understand if then $X$ and $Y$ {\it each} follow a Gamma law. This is related to Cram\'er's theorem which states that if $X$ and $Y$ are independent then $Z=X+Y$ follows a Gaussian law if and only if $X$ {\it and} $Y$ follow a Gaussian law. We prove that Cram\'er's theorem is true in the Gamma context for random variables leaving in a Wiener chaos of fixed order but the result is not true in general. We also give an asymptotic variant of our result.\\


\noindent
{\bf 2010 AMS Classification Numbers:} 60F05, 60H05, 91G70. }
\\

\section{Introduction}

\noindent Cram\'er's theorem (see \cite{HC}) says that the sum of two independent random variables is Gaussian if and only if each summand is Gaussian. One direction is elementary to prove, that is, given two independent random variables with Gaussian distribution, then their sum follows a Gaussian distribution. The second direction is less trivial and its proof requires powerful results from complex analysis (see \cite{HC}).
\\
\noindent In this paper, we treat the same problem for Gamma distributed random variables. A Gamma  random variable, denoted usually by $\Gamma (a, \lambda)$,  is a random variable with probability density function given by $f_{a,\lambda }(x)= \frac{\lambda ^{a}}{\Gamma (a)} x^{a-1} e^{-\lambda x} $ if $x>0$ and $f_{a, \lambda }(x)=0$ otherwise. The parameters $a$ and $\lambda $ are strictly positive and $\Gamma$ denotes the usual Gamma function.
\\
\noindent It is well known that if $X\sim \Gamma (a, \lambda) $ and $Y\sim \Gamma (b, \lambda) $  and $X$ is independent of $Y$, then $X+ Y$ follows the law $\Gamma (a+b, \lambda)$. The purpose of this paper is to understand the converse implication, i.e. whether or not (or under what conditions), if $X$ and $Y$ are two independent random variables such that $X+ Y\sim \Gamma (a+b, \lambda) $ and $\mathbf{E}(X)=\mathbf{E}\left(\Gamma (a, \lambda)\right), \mathbf{E}\left(X^{2}\right) =\mathbf{E}\left(\Gamma (a, \lambda ) ^{2}\right) $ and  $\mathbf{E}(Y)=\mathbf{E}\left(\Gamma (b, \lambda)\right), \mathbf{E}\left(Y^{2}\right) =\mathbf{E}\left(\Gamma (b, \lambda ) ^{2}\right) $, it holds that $X\sim \Gamma (a, \lambda ) $ and $Y\sim \Gamma (b, \lambda)$.
\\
\noindent We will actually focus our attention on the so-called centered Gamma distribution $F(\nu)$. We will call `centered Gamma' the random variables of the form
\begin{equation*}
F(\nu)\stackrel{\rm Law}{=}2G(\nu/2)-\nu, \quad \nu>0,
\end{equation*}
where $G(\nu/2):=F(\nu/2, 1)$ has a Gamma law with parameters $\nu/2 , 1$. This means that $\Gamma (\nu/2, 1)$ is a (a.s. strictly positive) random variable with density $g(x) = \frac{x^{\frac{\nu}{2}-1}{\rm e}^{-x}}{\Gamma(\nu/2)}\mathbf{1}_{(0,\infty)}(x)$. The characteristic function of the law $F(\nu)$ is given by
\begin{eqnarray}\label{fc}
\mathbf{E}\left(e^{i\lambda F(\nu)}\right) = \left(\frac{e^{-i\lambda}}{\sqrt{1-2i\lambda}}\right)^{\nu}, \mbox{\  \  \  \  } \lambda \in \mathbb{R}.
\end{eqnarray}
\\
\noindent We will find the following answer: if $X$ and $Y$ and two independent random variables, each leaving in a Wiener chaos of fixed order (and these orders are allowed to be different) then the fact that the sum $X+Y$ follows a centered Gamma distribution implies that $X$ and $Y$ each follow a Gamma distribution. On the other hand, for random variables having an infinite Wiener-It\^o chaos decomposition, the result is not true even in very particular cases (for so-called strongly independent random variables). We construct a counter-example to  illustrate this fact.
\\
\noindent Our tools are based on a criterium given in \cite{NoPe1} to characterize the random variables with Gamma distribution in terms of Malliavin calculus.
\\
\noindent Our paper is structured as follows. Section 2 contains some notations and preliminaries. In Section 3 we prove the Cram\'er theorem for Gamma distributed random variables in Wiener chaos of finite orders and we also give an asymptotic version of this result. In Section 4 we show that the result does not hold in the general case.

\section{Some notations and definitions}

Let $(W_{t})_{t\in T}$ be a classical Wiener process on a standard Wiener space $\left(\Omega,{\mathcal{F}},\mathbf{P}\right)$. If $f\in L^{2}(T^{n})$ with $n\geq 1$ integer, we introduce the multiple Wiener-It\^{o} integral of $f$ with respect to $W$. The basic references are the monographs \cite{Mal} or \cite{N}. Let $f\in {\mathcal{S}_{n}}$ be an elementary function with $n$
variables  that can be written as $
f=\sum_{i_{1},\ldots ,i_{n}}c_{i_{1},\ldots ,i_{n}}1_{A_{i_{1}}\times
\ldots \times A_{i_{n}}}$
where the coefficients satisfy $c_{i_{1},\ldots ,i_{n}}=0$ if two indices $i_{k}$ and $i_{l}$ are equal and the sets $A_{i}\in
{\mathcal{B}}(T)$ are pairwise disjoint. For  such a step function $f$ we define
\begin{equation*}
I_{n}(f)=\sum_{i_{1},\ldots ,i_{n}}c_{i_{1},\ldots ,i_{n}}W(A_{i_{1}})\ldots W(A_{i_{n}})
\end{equation*}
where we put $W(A)=\int_{0}^{1} 1_{A}(s)dW_{s}$. It can be seen that the application $
I_{n}$ constructed above from ${\mathcal{S}}_{n}$ to $L^{2}(\Omega )$ is an isometry on ${\mathcal{S}}_{n}$  in the sense
\begin{equation}
\mathbf{E}\left( I_{n}(f)I_{m}(g)\right) =n!\langle f,g\rangle_{L^{2}(T^{n})}\mbox{ if }m=n  \label{isom}
\end{equation}
and
\begin{equation*}
\mathbf{E}\left( I_{n}(f)I_{m}(g)\right) =0\mbox{ if }m\not=n.
\end{equation*}
\noindent Since the set ${\mathcal{S}_{n}}$ is dense in $L^{2}(T^{n})$ for every $n\geq 1$ the mapping $ I_{n}$ can be extended to an isometry from $L^{2}(T^{n})$ to $L^{2}(\Omega)$ and the above properties hold true for this extension.
\\
\noindent It also holds that $
I_{n}(f) = I_{n}\big( \tilde{f}\big)$ where $\tilde{f} $ denotes the symmetrization of $f$ defined by $$\tilde{f}(x_{1}, \ldots , x_{n}) =\frac{1}{n!} \sum_{\sigma}f(x_{\sigma (1) }, \ldots , x_{\sigma (n) } ),$$ $\sigma$ running over all permutations of $\left\{1,...,n\right\}$.
We will need the general formula for calculating products of Wiener chaos integrals of any orders $m,n$ for any symmetric integrands $f\in L^{2}(T^{m})$ and $g\in L^{2}(T^{n})$, which is
\begin{equation}
I_{m}(f)I_{n}(g)=\sum_{\ell=0}^{m\wedge n}\ell!\binom{m}{\ell}\binom{n}{\ell}%
I_{m+n-2\ell}(f\otimes_{\ell}g) \label{product}%
\end{equation}
where the contraction $f\otimes_{\ell}g$ is  defined by
\begin{eqnarray}
&&  (f\otimes_{\ell} g) ( s_{1}, \ldots, s_{m-\ell}, t_{1}, \ldots, t_{n-\ell
})\nonumber\\
&&  =\int_{T ^{m+n-2\ell} } f( s_{1}, \ldots, s_{m-\ell}, u_{1},
\ldots,u_{\ell})g(t_{1}, \ldots, t_{n-\ell},u_{1}, \ldots,u_{\ell})
du_{1}\ldots du_{\ell} . \label{contra}%
\end{eqnarray}
Note that the contraction $(f\otimes_{\ell} g) $ is an element of $L^{2}(T^{m+n-2\ell})$ but it is not necessarily symmetric. We will denote its symmetrization by $(f \tilde{\otimes}_{\ell} g)$.
\\
\noindent We recall that any square integrable random variable which is
measurable with respect to the $\sigma$-algebra generated by $W$ can
be expanded into an orthogonal sum of multiple stochastic integrals
\begin{equation}
\label{sum1} F=\sum_{n\geq0}I_{n}(f_{n})
\end{equation}
where $f_{n}\in L^{2}(T^{n})$ are (uniquely determined)
symmetric functions and $I_{0}(f_{0})=\mathbf{E}\left(F\right)  $.
\\
\noindent We denote by $D$  the Malliavin  derivative operator that acts on smooth functionals of the form $F=g(W(\varphi _{1}), \ldots , W(\varphi_{n}))$ (here $g$ is a smooth function with compact support and $\varphi_{i} \in L^{2}(T)$ for $i=1,..,n$)
\begin{equation*}
DF=\sum_{i=1}^{n}\frac{\partial g}{\partial x_{i}}(W(\varphi _{1}), \ldots , W(\varphi_{n}))\varphi_{i}.
\end{equation*}
We can define the $i$-th Malliavin derivative $D^{(i)}$ iteratively. The operator $D^{(i)}$ can be extended to the closure $\mathbb{D}^{p,2}$ of smooth functionals with respect to the norm
\begin{equation*}
\Vert F\Vert _{p,2}^{2} = \mathbf{E}\left(F^{2}\right)+ \sum_{i=1}^{p} \mathbf{E} \left(\Vert D^{i} F\Vert ^{2} _{L^{2}(T ^{i})}\right).
\end{equation*}
The adjoint of $D$ is denoted by $\delta $ and is called the divergence (or Skorohod) integral. Its domain $\mbox{Dom}(\delta)$  coincides with the class of stochastic processes $u\in L^{2}(\Omega \times T)$ such that
\begin{equation*}
\left| \mathbf{E}\left(\langle DF, u\rangle\right) \right| \leq c\Vert F\Vert _{2}
\end{equation*}
for all $F\in \mathbb{D}^{1,2}$ and $\delta (u)$ is the element of $L^{2}(\Omega)$ characterized by the duality relationship
\begin{equation*}
\mathbf{E}(F\delta (u))= \mathbf{E}\left(\langle DF, u\rangle\right).
\end{equation*}
For adapted integrands, the divergence integral coincides with the classical It\^o integral.
\\
\noindent Let $L$ be the Ornstein-Uhlenbeck operator defined on $\mbox{Dom}(L)= \mathbb{D}^{2,2}$. We have
\begin{equation*}
LF=-\sum_{n\geq 0} nI_{n}(f_{n})
\end{equation*}
if $F$ is given by (\ref{sum1}). There exists a connection between $\delta, D $ and $L$ in the sense that a random variable $F$ belongs to the domain of $L$ if and only if $F\in \mathbb{D}^{1,2}$ and $DF \in \mbox{Dom}(\delta)$ and then $\delta DF=-LF$.
\noindent Let us consider a multiple stochastic integral $I_{q}(f)$ with symmetric kernel $f\in L^{2}(T^{q})$. We denote the Malliavin derivative of $I_{q}(f)$ by $DI_{q}(f)$. We have $$D_{\theta}I_{q}(f) = qI_{q-1}(f^{(\theta)}),$$ where $f^{(\theta)} = f(t_{1},...,t_{q-1},\theta)$ is the $(q-1)^{\mbox{\tiny{th}}}$ order kernel obtained by parametrizing the $q^{\mbox{\tiny{th}}}$ order kernel $f$ by one of the variables.
\\
\noindent For any random variable $X,Y \in \mathbb{D}^{1,2}$ we use the following notations
\begin{equation*}\label{gx}
G_{X} = \langle DX, -DL ^{-1}X\rangle_{L^{2}(T)}
\end{equation*}
and
\begin{equation*}\label{gxy}
G_{X,Y} = \langle DX, -DL ^{-1}Y\rangle_{L^{2}(T)}.
\end{equation*}
\\
\noindent The following facts are key points in our proofs:
\begin{description}
\item{ {\bf Fact 1: }} Let $X=I_{q_{1}}(f)  $ and $Y= I_{q_{2}} (g)$ where $f\in L^{2}(T^{q_{1}})$ and $g\in L^{2}(T^{q_{2}} )$ are symmetric functions. Then $X$ and $Y$ are independent if and only if (see \cite{UsZa})
    \begin{equation*}
    f\otimes _{1} g =0 \mbox{ a.e. on } T^{q_{1}+q_{2}-2 } .
    \end{equation*}

    \item{{\bf Fact 2: } } Let $X=I_{q}(f)  $ with  $f\in L^{2}(T^{q})$ symmetric. Assume that $\mathbf{E}\left(X^{2}\right)=\mathbf{E}(F(\nu ) ^{2}) =2\nu$.  Then $X$ follows a centered Gamma law $F(\nu)$ with $\nu >0$ if and only if (see \cite{NoPe2})
        \begin{equation*}
        \Vert DX \Vert ^{2} _{L^{2}(T) } -2q X-2q \nu =0 \mbox{ almost surely. }
        \end{equation*}
\item{{\bf Fact 3: }} Let $(f_{k}) _{k\geq 1}$ be a sequence in $L^{2}(T^{q} ) $ such that $\mathbf{E}\left(I_{q} (f_{k}) ^{2}\right)\underset{k \rightarrow +\infty}{\longrightarrow} 2\nu. $ Then the sequence $X_{k}=  I_{q}(f_{k}) $ converges in  distribution, as $k\to \infty$, to a Gamma law, if and only if (see \cite{NoPe2})
    \begin{equation*}
     \Vert DX _{k}\Vert ^{2} _{L^{2}(T) } -2q X_{k} -2q \nu \underset{k \rightarrow +\infty}{\longrightarrow} 0 \mbox{ in } L^{2}(\Omega).
    \end{equation*}

\end{description}

\begin{Remark}
In this particular paper, we will restrict ourselves to an underlying Hilbert space (to the Wiener process we will be working with in the upcoming sections) of the form $\HH = L^{2}(T)$ for the sake of simplicity. However, all the results presented in the upcoming sections remain valid on a more general separable Hilbert space as the underlying space.
\end{Remark}

\section{(Asymptotic) Cram\'er theorem for multiple integrals}

In this section, we will prove Cram\'er's theorem for random variables living in fixed Wiener chaoses. More precisely, our context is as follows: we assume that $X= I_{q_{1}}(f)$ and $Y=I_{q_{2}} (h)$ and $X, Y$ are independent. We also assume that $\mathbf{E}\left(X^{2}\right)= \mathbf{E}\left(F(\nu_{1}) ^{2}\right) =2\nu _{1} $ and $\mathbf{E}\left(Y^{2}\right)= \mathbf{E}\left(F(\nu _{2} ^{2})\right) =2\nu _{2}$. Here $\nu, \nu_{1}, \nu_{2} $ denotes three strictly positive numbers such that $\nu_{1}+ \nu_{2}=\nu$. We assume that $X+Y$ follows a Gamma law $F(\nu)$ and we will prove that $X\sim F(\nu _{1}) $ and $Y\sim F(\nu _{2})$.
\\
\noindent Let us first give the two following auxiliary lemmas that will be useful throughout the paper.
\begin{lemma}
\label{lemmeIndep}
Let $q_{1},q_{2} \geq 1$ be integers, and let $X = I_{q_{1}}(f)$ and $Y = I_{q_{2}}(h)$, where $f \in L^{2}(T^{ q_{1}})$ and $h \in L^{2}(T^{ q_{2}})$ are symmetric functions. Assume moreover that $X$ and $Y$ are independent. Then, we have $DX \bot DY$, $X \bot DY$ and $Y \bot DX$.
\end{lemma}
{\bf Proof: }
From Fact 1 in Section 2,  $f \otimes_{1} h = 0$ a.e on $T^{q_{1}+q_{2}-2}$  and by extension $f \otimes_{r} h = 0$ a.e  on $T^{q_{1}+q_{2}-2r}$ for every  $1\leq r \leq q_{1} \wedge q_{2}$.
We will now prove that for every $\theta, \psi \in T$, we also have $f^{(\theta)} \otimes_{1} h^{(\psi)} = 0$ a.e on $T^{q_{1}+q_{2} -4}$, $f^{(\theta)} \otimes_{1} h = 0$ a.e on $T^{q_{1}+q_{2}-3}$ and $f \otimes_{1} h^{(\psi)} = 0$ a.e. on $T^{q_{1}+q_{2}-3}$. Indeed, we have
\begin{eqnarray*}
\left( f^{(\theta)} \otimes_{1} h^{(\psi)} \right) (t_{1}, \ldots,t_{q_{1}-2}, s_{1}, \ldots , s_{q_{2}-2})&=& \int_{T}f(t_{1},...,t_{q_{1}-2},u,\theta)h(s_{1},...,s_{q_{2}-2},u,\psi)du
\\
&=& 0
\end{eqnarray*}
as a particular case of $f \otimes_{1} h = 0$ a.e.. By extension, we also have $f^{(\theta)} \otimes_{r} h^{(\psi)} = 0$ for $1\leq r \leq (q_{1}-1) \wedge (q_{2}-1)$. Similarly,
\begin{eqnarray}
\left( f^{(\theta)} \otimes_{1} h\right) (t_{1}, \ldots,t_{q_{1}-2}, s_{1}, \ldots , s_{q_{2}-1}) &=& \int_{T}f(t_{1},...,t_{q_{1}-2},u,\theta)h(s_{1}, \ldots , s_{q_{2}-1},u)du.
\nonumber \\
&=& 0 \label{b1}
\end{eqnarray}
 and clearly $f^{(\theta)} \otimes_{r} h = 0$ for $1\leq r \leq (q_{1}-1) \wedge q_{2}$.
Given the symmetric roles played by $f$ and $h$, we also have $f \otimes_{1} h^{(\psi)} = 0$ and then  $f \otimes_{r} h^{(\psi)} = 0$ for $1\leq r \leq q_{1} \wedge (q_{2}-1)$.
\\
\noindent Let us now prove that $DX \bot DY$. Since for every $\theta, \psi \in T$, $D_{\theta} X= q_{1} I_{q_{1}}(f^{(\theta)})$ and $D_{\psi} Y= q_{2} I_{q_{2}-1}(h^{(\psi )})$, it suffices to
show that the random variables $I_{q_{1}}(f^{(\theta)})$ and $I_{q_{2}-1}(h^{(\psi )})$ are independent. To do this, we will use the criterium for the independence of multiple integrals given in \cite{UsZa}. We need to check that
$f^{(\theta ) } \otimes _{1} h^{(\psi)} =0$ a.e. on $T^{q_{1}+q_{2}-4}$ and this follows from above.
\\
\noindent It remains to prove that $X \bot DY$ and $DX \bot Y$. Given the symmetric roles played by $X$ and $Y$, we will only prove that $X \bot DY$. That is equivalent to the independence of the random variables $I_{q_{1}-1}(f^{(\theta)})$ and $I_{q_{2}}(h)$  for every $\theta \in T$, which follows from \cite{UsZa} (see Fact 1 in Section 2) and  (\ref{b1}). Thus, we have $X \bot DY$ and $DX \bot Y$.
\qed

\vskip0.2cm

\noindent Let us recall the following definition (see \cite{Tu}).
\begin{definition} Two random variables $X=\sum_{n\geq 0} I_{n}(f_{n})$
and $Y=\sum_{m\geq 0} I_{m}(h_{m}) $ are called {\it strongly independent} if  for every $m,n\geq 0$,
 the random variables $I_{n}(f_{n})$ and $I_{m}(h_{m})$ are  independent.
\end{definition}
We have the following lemma about strongly independent random variables.
\begin{lemma}
\label{lemmeG}
Let $X=\sum_{n\geq 0} I_{n}(f_{n})$ and $Y=\sum_{m\geq 0} I_{m}(h_{m})$ ($f_{n}\in L^{2}(T^{n}), h_{m}\in L^{2}(T^{m})$ symmetric for every $n,m\geq 1$) be two centered random variables in the space $\mathbb{D}^{1,2}$. Then, if $X$ and $Y$ are strongly independent, we have $$\langle DX, -DL ^{-1}Y\rangle_{L^{2}(T)} = \langle DY, -DL ^{-1}X\rangle_{L^{2}(T)} = 0.$$
\end{lemma}
{\bf Proof: }
We have, for every $\theta \in T$, $$D_{\theta}X = \sum_{n\geq 1}n I_{n-1}(f_{n}^{(\theta)}) \mbox{ and } -D_{\theta}L ^{-1}Y = \sum_{m\geq 1} I_{m-1}(h_{m}^{(\theta)}).$$ Therefore, we can write
\begin{eqnarray*}
\langle DX, -DL ^{-1}Y\rangle_{L^{2}(T)} &=& \sum_{n,m\geq 1}n\int_{T} I_{n-1}(f_{n}(t_{1},...,t_{n-1},\theta))I_{m-1}(h_{m}(t_{1},...,t_{m-1},\theta))d\theta
\\
&=& \sum_{n,m\geq 1}n\int_{T}\sum_{r=0}^{(n-1) \wedge (m-1)}r!\binom{n-1}{r}\binom{m-1}{r}I_{n+m-2r-2}(f_{n}^{(\theta)} \otimes_{r} h_{m}^{(\theta)})d\theta.
\end{eqnarray*}
The strong independence of $X$ and $Y$ gives us that $f_{n}^{(\theta)} \otimes_{r} h_{m}^{(\theta)} = 0$ for every $1 \leq r \leq (n-1) \wedge (m-1)$. Thus, we obtain
\begin{eqnarray*}
\langle DX, -DL ^{-1}Y\rangle_{L^{2}(T)} = \sum_{n,m\geq 1}n\int_{T}I_{n+m-2}(f_{n}^{(\theta)} \otimes h_{m}^{(\theta)})d\theta.
\end{eqnarray*}
Using a Fubini type result, we can write
\begin{eqnarray*}
\langle DX, -DL ^{-1}Y\rangle_{L^{2}(T)} &=& \sum_{n,m\geq 1}nI_{n+m-2}(\int_{T}f_{n}^{(\theta)} \otimes h_{m}^{(\theta)}d\theta)
\\
&=& \sum_{n,m\geq 1}nI_{n+m-2}(f_{n} \otimes_{1} h_{m}).
\end{eqnarray*}
Again, the strong independence of $X$ and $Y$ gives us that $f_{n} \otimes_{1} h_{m} = 0$ a.e and we finally obtain $
\langle DX, -DL ^{-1}Y\rangle_{L^{2}(T)} = 0,
$ and similarly $\langle DY, -DL ^{-1}X\rangle_{L^{2}(T)} = 0.$ \qed

\vskip0.3cm

\noindent Let us first remark that that the Cram\'er theorem holds for random variables in the same   Wiener  chaos of fixed order.

\begin{prop}
Let $X=I_{m}(f)$ and $Y=I_{m}(h)$ with $m\geq 2$  fixed and $f,h$ symmetric functions in $L^{2}(T^{m})$. Then $X+Y= I_{m}(f+h)$.  Fix $\nu _{1}, \nu _{2}, \nu >0$ such that $\nu_{1}+ \nu _{2}= \nu$.  Assume that $X+Y$ follows the law $F(\nu)$ and $X$ is independent of $Y$. Also suppose that $\mathbf{E} \left(X^{2}\right)=\mathbf{E}\left(F(\nu_{1} )^{2}\right)= 2\nu_{1} $ and $\mathbf{E}\left(Y^{2}\right)=\mathbf{E}\left(F(\nu_{2})^{2}\right) = 2\nu _{2}$.
 Then $X\sim F\left(\nu _{1}\right)$  and  $Y\sim F\left(\nu _{2}\right)$.
 \end{prop}
{\bf Proof: } By a result in \cite{NoPe2} (see Fact 2 in Section 2),  $X+Y$ follows the law $F(\nu)$  is equivalent to
\begin{equation*}
\left|\left|D I_{m}(f+h)\right| \right| ^{2} _{L^{2}(T)}-2m I_{m}(f+h) -2m \nu = 0 \mbox{ a.s. }.
\end{equation*}
On the other hand
\begin{eqnarray*}
&&\mathbf{E}\left( \left|\left|D I_{m}(f+h)\right| \right| ^{2} _{L^{2}(T)}-2m I_{m}(f+h) -2m \nu\right) ^{2} \\
&=&\mathbf{E}\left(\left(   \left|\left|D I_{m}(f)\right| \right| ^{2} _{L^{2}(T)} + \left|\left|D I_{m}(h)\right| \right| ^{2} _{L^{2}(T)}
+2\langle DI_{m}(f), DI_{m}(h)\rangle _{L^{2}(T)} \right.\right. \\
&& \left.\left. -2mI_{m}(f)-2mI_{m}(h) -2m(\nu_{1}+ \nu_{2} ) \right)^{2}\right)  \\
&=& \mathbf{E}\left(\left( \left|\left|D I_{m}(f)\right| \right| ^{2} _{L^{2}(T)} -2m I_{m}(f) -2m\nu _{1} \right) ^{2}\right) + \mathbf{E}\left(\left( \left|\left|D I_{m}(h)\right| \right| ^{2} _{L^{2}(T)} -2m I_{m}(h) -2m\nu _{2} \right) ^{2}\right) \\
&&+ \mathbf{E}\left( \left( \left|\left|D I_{m}(f)\right| \right| ^{2} _{L^{2}(T)} -2m I_{m}(f) -2m\nu _{1} \right) \left( \left|\left|D I_{m}(h)\right| \right| ^{2} _{L^{2}(T)} -2m I_{m}(h) -2m\nu _{2} \right)\right).
\end{eqnarray*}
Above we used the fact that $\langle DI_{m}(f), DI_{m}(h)\rangle _{L^{2}(T)}=0 $ as a consequence of Lemma \ref{lemmeIndep}. It is also easy to remark that, from Lemma \ref{lemmeIndep}
\begin{eqnarray*}
&&\mathbf{E}\left( \left( \left|\left|D I_{m}(f)\right| \right| ^{2} _{L^{2}(T)} -2m I_{m}(f) -2m\nu _{1} \right) \left( \left|\left|D I_{m}(h)\right| \right| ^{2} _{L^{2}(T)} -2m I_{m}(h) -2m\nu _{2} \right)\right)\\
&=& \mathbf{E} \left( \left|\left|D I_{m}(f)\right| \right| ^{2} _{L^{2}(T)} -2m I_{m}(f) -2m\nu _{1} \right)\mathbf{E}\left( \left|\left|D I_{m}(h)\right| \right| ^{2} _{L^{2}(T)} -2m I_{m}(h) -2m\nu _{2} \right)=0.
\end{eqnarray*}
We will obtain that
\begin{equation*}
\mathbf{E}\left(\left( \left|\left|D I_{m}(f)\right| \right| ^{2} _{L^{2}(T)} -2m I_{m}(f) -2m\nu _{1} \right) ^{2}\right) =\mathbf{E}\left(\left( \left|\left|D I_{m}(h)\right| \right| ^{2} _{L^{2}(T)} -2m I_{m}(h) -2m\nu _{2} \right) ^{2}\right)=0
\end{equation*}
and consequently $X\sim F(\nu _{1}) $ and $Y \sim F(\nu _{2})$. \qed

\begin{remark}
Using Fact 3 in Section 2, an asymptotic variant of the above result can be stated. We will state it here because it is a particular case of Theorem \ref{cramerAsym} proved later in our paper.
\end{remark}

\noindent Theorem 1.2 in \cite{NoPe2} gives a characterization of (asymptotically) centered Gamma random variable which are given by a multiple Wiener-It\^o integral. There is not such a characterization  for random variable leaving in  a finite or infinite sum of Wiener chaos; only an upper bound for the distance between the law of a random variable in $\mathbb{D} ^{1,2}$ and the Gamma distribution has been proven in \cite{NoPe1}, Theorem 3.11. It turns out, that for the case of a sum of independent multiple integrals, it is possible to characterize the relation between its distribution and the Gamma distribution. We will prove this fact in the following theorem.

\begin{theorem}
\label{theoGeneral}
\noindent Fix $\nu_{1},\nu_{2},\nu > 0$ such that $\nu_{1}+\nu_{2}=\nu$ and let $F(\nu)$ be a real-valued random variable with characteristic function given by (\ref{fc}).
\noindent Fix two even integers $q_{1} \geq 2$ and $q_{2} \geq 2$. For any symmetric kernels $f \in L^{2}(T^{ q_{1}})$ and $h \in L^{2}(T^{ q_{2}})$ such that
\begin{eqnarray}
\label{momentordre2}
\mathbf{E}\left(I_{q_{1}}(f)^{2}\right) = q_{1}!\left\|f\right\|^2_{L^{2}(T^{q_{1}})} = 2\nu_{1} \mbox{\  \  \  and\  \  \  } \mathbf{E}\left(I_{q_{2}}(h)^{2}\right) = q_{2}!\left\|h\right\|^2_{L^{2}(T^{q_{2}})} = 2\nu_{2},
\end{eqnarray}
and such that $X = I_{q_{1}}(f)$ and $Y = I_{q_{2}}(h)$ are independent, define the random variable
\begin{eqnarray*}
Z = X + Y = I_{q_{1}}(f) + I_{q_{2}}(h).
\end{eqnarray*}
Under those conditions, the following two conditions are equivalent:
\begin{enumerate}[(i)]
\item $\mathbf{E}\left(\left(2\nu + 2Z - \left\langle DZ,-DL^{-1}Z\right\rangle_{L^{2}(T)}\right)^{2}\right) = 0$, where $D$ is the Malliavin derivative operator and $L$ is the infinitesimal generator of the Ornstein-Uhlenbeck semigroup;
\item $Z \stackrel{\rm Law}{=} F(\nu)$;
\end{enumerate}
\end{theorem}
{\bf Proof: }
\textit{ Proof of $(ii) \rightarrow (i)$}. Suppose that $Z\sim F(\nu)$. We easily obtain that
\begin{equation}\label{equiv2}
\mathbf{E}\left(Z^{3}\right) = \mathbf{E}\left(F(\nu)^{3}\right) = 8\nu \mbox{\   \    and \   \   }\mathbf{E}\left(Z^{4}\right) = \mathbf{E}\left(F(\nu)^{4}\right) = 12\nu^{2} + 48\nu.
\end{equation}
Consequently,
\begin{equation}\label{equiv1}
\mathbf{E}\left(Z^{4}\right)  -12\mathbf{E}\left(Z^{3}\right) = \mathbf{E} \left(F(\nu) ^{4}\right) -12\mathbf{E}\left(F(\nu ) ^{3}\right)=12\nu ^{2} -48 \nu.
\end{equation}
Then we  will  use the fact that for every multiple integral $I_{q}(f)$
\begin{eqnarray}
\label{cubeIntMul}
\mathbf{E}\left(I_{q}(f)^{3}\right) = q!(q/2)!\binom{q}{q/2}^{2}\left\langle f, f \widetilde{\otimes}_{q/2} f\right\rangle_{L^{2}(T^{q})}.
\end{eqnarray}
and
\begin{eqnarray}
\label{puiss4IntMul2}
\mathbf{E}\left(I_{q}(f)^{4}\right) = 3\left[q!\left\|f\right\|_{L^{2}(T^{q})}^{2}\right]^{2}+\frac{3}{q}\sum_{p=1}^{q-1}q^{2}(p-1)!\binom{q-1}{p-1}^{2}p!\binom{q}{p}^{2}(2q-2p)!\left\|f \widetilde{\otimes}_{p} f\right\|_{L^{2}(T^{2(q - p)})}^{2}.
\end{eqnarray}
\noindent We will now compute $\mathbf{E}\left(Z^{3}\right)$, $\mathbf{E}\left(Z^{4}\right)$ and $\mathbf{E}\left(Z^{4}\right) - 12\mathbf{E}\left(Z^{3}\right)$ by using the above two relations (\ref{cubeIntMul}) and (\ref{puiss4IntMul2}).
We have $Z^{2} = (I_{q_{1}}(f) + I_{q_{2}}(h))^{2} = I_{q_{1}}(f)^{2} + I_{q_{2}}(h)^{2} +2 I_{q_{1}}(f)I_{q_{2}}(h)$ and thus, by using the independence between $I_{q_{1}}(f)$ and $I_{q_{2}}(h)$,
\begin{eqnarray*}
\mathbf{E}\left(Z^{3}\right) &=& \mathbf{E}\left(I_{q_{1}}(f)^{3}\right) + \mathbf{E}\left(I_{q_{2}}(h)^{3}\right).
\end{eqnarray*}
Using relation (\ref{cubeIntMul}), we can write
\begin{eqnarray}
\label{Z3}
\mathbf{E}\left(Z^{3}\right) = q_{1}!(q_{1}/2)!\binom{q_{1}}{q_{1}/2}^{2}\left\langle f, f \widetilde{\otimes}_{q_{1}/2} f\right\rangle_{L^{2}(T^{ q_{1}})} + q_{2}!(q_{2}/2)!\binom{q_{2}}{q_{2}/2}^{2}\left\langle h, h \widetilde{\otimes}_{q_{2}/2} h\right\rangle_{L^{2}(T^{ q_{2}})}.
\end{eqnarray}
For $\mathbf{E}\left(Z^{4}\right)$, we combine relations (\ref{momentordre2}) and (\ref{puiss4IntMul2}) with the independence between $I_{q_{1}}(f)$ and $I_{q_{2}}(h)$ to obtain
\begin{eqnarray*}
\mathbf{E}\left(Z^{4}\right) &=& \mathbf{E}\left(Z^{2}Z^{2}\right) = \mathbf{E}\left(I_{q_{1}}(f)^{4}\right) + \mathbf{E}\left(I_{q_{2}}(h)^{4}\right) + 6\mathbf{E}\left(I_{q_{1}}(f)^{2}I_{q_{2}}(h)^{2}\right)
\\
&=& 3\left[q_{1}!\left\|f\right\|_{L^{2}(T^{q_{1}})}^{2}\right]^{2}+\frac{3}{q_{1}}\sum_{p=1}^{q_{1}-1}q_{1}^{2}(p-1)!\binom{q_{1}-1}{p-1}^{2}p!\binom{q_{1}}{p}^{2}(2q_{1}-2p)!\left\|f \widetilde{\otimes}_{p} f\right\|_{L^{2}(T^{ 2(q_{1} - p)})}^{2}
\\
&& + 3\left[q_{2}!\left\|h\right\|_{L^{2}(T^{q_{2}})}^{2}\right]^{2}+\frac{3}{q_{2}}\sum_{p=1}^{q_{2}-1}q_{2}^{2}(p-1)!\binom{q_{2}-1}{p-1}^{2}p!\binom{q_{2}}{p}^{2}(2q_{2}-2p)!\left\|h \widetilde{\otimes}_{p} h\right\|_{L^{2}(T^{ 2(q_{2} - p)})}^{2}
\\
&& + 24\nu_{1}\nu_{2}.
\end{eqnarray*}
Using the fact that $q_{1}!\left\|f\right\|_{L^{2}(T^{ q_{1}})}^{2} = 2\nu_{1}$ and $q_{2}!\left\|h\right\|_{L^{2}(T^{ q_{2}})}^{2} = 2\nu_{2}$, we can write
\begin{eqnarray}
\label{Z412Z3}
\mathbf{E}\left(Z^{4}\right) - 12\mathbf{E}\left(Z^{3}\right) &=& 12\nu_{1}^{2}+12\nu_{2}^{2} - 48\nu_{1} - 48\nu_{2} + 24\nu_{1}\nu_{2}\nonumber
\\
&& +\frac{3}{q_{1}}\sum_{p=1, p \neq q_{1}/2}^{q_{1}-1}q_{1}^{2}(p-1)!\binom{q_{1}-1}{p-1}^{2}p!\binom{q_{1}}{p}^{2}(2q_{1}-2p)!\left\|f \widetilde{\otimes}_{p} f\right\|_{L^{2}(T^{ 2(q_{1} - p)})}^{2}\nonumber
\\
&& +\frac{3}{q_{2}}\sum_{p=1, p \neq q_{2}/2}^{q_{2}-1}q_{2}^{2}(p-1)!\binom{q_{2}-1}{p-1}^{2}p!\binom{q_{2}}{p}^{2}(2q_{2}-2p)!\left\|h \widetilde{\otimes}_{p} h\right\|_{L^{2}(T^{ 2(q_{2} - p)})}^{2}\nonumber
\\
&& + 24q_{1}!\left\|f\right\|_{L^{2}(T^{ q_{1}})}^{2} + 3q_{1}(q_{1}/2 - 1)!\binom{q_{1}-1}{q_{1}/2-1}^{2}(q_{1}/2)!\binom{q_{1}}{q_{1}/2}^{2}q_{1}!\left\|f \widetilde{\otimes}_{q_{1}/2} f\right\|_{L^{2}(T^{ q_{1}})}^{2}\nonumber
\\
&& + 24q_{2}!\left\|h\right\|_{L^{2}(T^{ q_{2}})}^{2} + 3q_{2}(q_{2}/2 - 1)!\binom{q_{2}-1}{q_{2}/2-1}^{2}(q_{2}/2)!\binom{q_{2}}{q_{2}/2}^{2}q_{2}!\left\|h \widetilde{\otimes}_{q_{2}/2} h\right\|_{L^{2}(T^{ q_{2}})}^{2}\nonumber
\\
&& - 12q_{1}!(q_{1}/2)!\binom{q_{1}}{q_{1}/2}^{2}\left\langle f, f \widetilde{\otimes}_{q_{1}/2} f\right\rangle_{L^{2}(T^{ q_{1}})} \nonumber
\\
&& - 12q_{2}!(q_{2}/2)!\binom{q_{2}}{q_{2}/2}^{2}\left\langle h, h \widetilde{\otimes}_{q_{2}/2} h\right\rangle_{L^{2}(T^{ q_{2}})}.
\end{eqnarray}
Recall that $\nu = \nu_{1} + \nu_{2}$ and note that $12\nu_{1}^{2}+12\nu_{2}^{2} - 48\nu_{1} - 48\nu_{2} + 24\nu_{1}\nu_{2} = 12\nu^{2} - 48\nu$. Also note that
\begin{eqnarray*}
&& 24q_{1}!\left\|f\right\|_{L^{2}(T^{ q_{1}})}^{2} + 3q_{1}(q_{1}/2 - 1)!\binom{q_{1}-1}{q_{1}/2-1}^{2}(q_{1}/2)!\binom{q_{1}}{q_{1}/2}^{2}q_{1}!\left\|f \widetilde{\otimes}_{q_{1}/2} f\right\|_{L^{2}(T^{ q_{1}})}^{2}
\\
&& - 12q_{1}!(q_{1}/2)!\binom{q_{1}}{q_{1}/2}^{2}\left\langle f, f \widetilde{\otimes}_{q_{1}/2} f\right\rangle_{L^{2}(T^{ q_{1}})}
\\
&=& \frac{3}{2}\frac{(q_{1}!)^{5}}{\left((q_{1}/2)!\right)^{6}}\left\|f \widetilde{\otimes}_{q_{1}/2} f - c_{q_{1}}f\right\|_{L^{2}(T^{ q_{1}})}^{2},
\end{eqnarray*}
where $c_{q_{1}}$ is defined by  $c_{q_{1}} = \frac{1}{(q_{1}/2)!\binom{q_{1} - 1}{q_{1}/2 -1}^{2}} = \frac{4}{(q_{1}/2)!\binom{q_{1}}{q_{1}/2}^{2}}
$
 and a similar relation holds for the function $h$ with $q_{2}, c_{q_{2}}$ instead of $q_{1}, c_{q_{1}}$ respectively, where $
 c_{q_{2}} = \frac{1}{(q_{2}/2)!\binom{q_{2} - 1}{q_{2}/2 -1}^{2}} = \frac{4}{(q_{2}/2)!\binom{q_{2}}{q_{2}/2}^{2}}.
$
\begin{eqnarray*}
\mathbf{E}\left(Z^{4}\right) - 12\mathbf{E}\left(Z^{3}\right) &=& 12\nu^{2} - 48\nu
\\
&& +\frac{3}{q_{1}}\sum_{p=1, p \neq q_{1}/2}^{q_{1}-1}q_{1}^{2}(p-1)!\binom{q_{1}-1}{p-1}^{2}p!\binom{q_{1}}{p}^{2}(2q_{1}-2p)!\left\|f \widetilde{\otimes}_{p} f\right\|_{L^{2}(T^{ 2(q_{1} - p)})}^{2}
\\
&& + \frac{3}{2}\frac{(q_{1}!)^{5}}{\left((q_{1}/2)!\right)^{6}}\left\|f \widetilde{\otimes}_{q_{1}/2} f - c_{q_{1}}f\right\|_{L^{2}(T^{ q_{1}})}^{2}
\\
&& +\frac{3}{q_{2}}\sum_{p=1, p \neq q_{2}/2}^{q_{2}-1}q_{2}^{2}(p-1)!\binom{q_{2}-1}{p-1}^{2}p!\binom{q_{2}}{p}^{2}(2q_{2}-2p)!\left\|h \widetilde{\otimes}_{p} h\right\|_{L^{2}(T^{ 2(q_{2} - p)})}^{2}
\\
&& + \frac{3}{2}\frac{(q_{2}!)^{5}}{\left((q_{2}/2)!\right)^{6}}\left\|h \widetilde{\otimes}_{q_{2}/2} h - c_{q_{2}}h\right\|_{L^{2}(T^{ q_{2}})}^{2}.
\end{eqnarray*}
From $(ii)$, it follows that
\begin{eqnarray*}
&& \frac{3}{q_{1}}\sum_{p=1, p \neq q_{1}/2}^{q_{1}-1}q_{1}^{2}(p-1)!\binom{q_{1}-1}{p-1}^{2}p!\binom{q_{1}}{p}^{2}(2q_{1}-2p)!\left\|f \widetilde{\otimes}_{p} f\right\|_{L^{2}(T^{ 2(q_{1} - p)})}^{2}
\\
&& + \frac{3}{2}\frac{(q_{1}!)^{5}}{\left((q_{1}/2)!\right)^{6}}\left\|f \widetilde{\otimes}_{q_{1}/2} f - c_{q_{1}}f\right\|_{L^{2}(T^{ q_{1}})}^{2}
\\
&& +\frac{3}{q_{2}}\sum_{p=1, p \neq q_{2}/2}^{q_{2}-1}q_{2}^{2}(p-1)!\binom{q_{2}-1}{p-1}^{2}p!\binom{q_{2}}{p}^{2}(2q_{2}-2p)!\left\|h \widetilde{\otimes}_{p} h\right\|_{L^{2}(T^{ 2(q_{2} - p)})}^{2}
\\
&& + \frac{3}{2}\frac{(q_{2}!)^{5}}{\left((q_{2}/2)!\right)^{6}}\left\|h \widetilde{\otimes}_{q_{2}/2} h - c_{q_{2}}h\right\|_{L^{2}(T^{ q_{2}})}^{2}
= 0,
\end{eqnarray*}
which leads to the conclusion as all the summands are positive, that is
\begin{eqnarray}
&&\left\|f \widetilde{\otimes}_{q_{1}/2} f - c_{q_{1}}f\right\|_{L^{2}(T^{ q_{1}})} = \left\|h \widetilde{\otimes}_{q_{2}/2} h - c_{q_{2}}h\right\|_{L^{2}(T^{q_{2}})} = 0 \mbox{ and } \nonumber \\
&& \left\|f \widetilde{\otimes}_{p} f\right\|_{L^{2}(T^{ 2(q_{1} - p)})} = \left\|h \widetilde{\otimes}_{r} h\right\|_{L^{2}(T^{ 2(q_{2} - p)})} = 0 \label{equiv3}
 \end{eqnarray}
 for every $p = 1,...,q_{1} - 1$ such that $p \neq q_{1}/2$ and for every $r = 1,...,q_{2} - 1$ such that $r \neq q_{2}/2$; This implies
\begin{eqnarray}
&&\left\|f \widetilde{\otimes}_{q_{1}/2} f - c_{q_{1}}f\right\|_{L^{2}(T^{q_{1}})} = \left\|h \widetilde{\otimes}_{q_{2}/2} h - c_{q_{2}}h\right\|_{L^{2}(T^{q_{2}})} = 0 \mbox{ and  }\nonumber \\
&& \left\|f \otimes_{p} f\right\|_{L^{2}(T^{ 2(q_{1} - p)})} = \left\|h \otimes_{r} h\right\|_{L^{2}(T^{ 2(q_{2} - p)})} = 0 \label{equiv4}
 \end{eqnarray}
 for every $p = 1,...,q_{1} - 1$ such that $p \neq q_{1}/2$ and for every $r = 1,...,q_{2} - 1$ such that $r \neq q_{2}/2$ (see \cite{NoPe2}, Theorem 1.2.).
\\
\noindent We will compute $\mathbf{E}\left(\left(2\nu + 2Z - G_{Z}\right)^{2}\right)$. Let us start with $G_{Z}$.
\begin{eqnarray*}
G_{Z} &=& \left\langle DZ, -DL^{-1}Z\right\rangle_{L^{2}(T)} = \left\langle DI_{q_{1}}(f) + DI_{q_{2}}(h), -DL^{-1}I_{q_{1}}(f) -DL^{-1}I_{q_{2}}(h)\right\rangle_{L^{2}(T)}
\\
&=& \left\langle DI_{q_{1}}(f) , -DL^{-1}I_{q_{1}}(f)\right\rangle_{L^{2}(T)} + \left\langle DI_{q_{2}}(h) , -DL^{-1}I_{q_{2}}(h)\right\rangle_{L^{2}(T)}
\\
& +& \left\langle DI_{q_{1}}(f) , -DL^{-1}I_{q_{2}}(h)\right\rangle_{L^{2}(T)} + \left\langle DI_{q_{2}}(h) , -DL^{-1}I_{q_{1}}(f)\right\rangle_{L^{2}(T)}.
\end{eqnarray*}
From Lemma \ref{lemmeG}, it follows that $\left\langle DI_{q_{1}}(f) , -DL^{-1}I_{q_{2}}(h)\right\rangle_{L^{2}(T)} = \left\langle DI_{q_{2}}(h) , -DL^{-1}I_{q_{1}}(f)\right\rangle_{L^{2}(T)}=0$. Thus,
\begin{eqnarray*}
G_{Z} =q_{1}^{-1}\left\|DI_{q_{1}}(f)\right\|_{L^{2}(T)}^{2} + q_{2}^{-1}\left\|DI_{q_{2}}(h)\right\|_{L^{2}(T)}^{2}.
\end{eqnarray*}
It follows that
\begin{eqnarray*}
&& \mathbf{E}\left(\left(2\nu + 2Z - G_{Z}\right)^{2}\right)
\\
&& = \mathbf{E}\left(\left(2\nu_{1} + 2\nu_{2} + 2I_{q_{1}}(f) + 2I_{q_{2}}(h) - q_{1}^{-1}\left\|DI_{q_{1}}(f)\right\|_{L^{2}(T)}^{2} - q_{2}^{-1}\left\|DI_{q_{2}}(h)\right\|_{L^{2}(T)}^{2}\right)^{2}\right)
\\
&& = \mathbf{E}\left(\left(q_{1}^{-1}\left\|DI_{q_{1}}(f)\right\|_{L^{2}(T)}^{2} - 2I_{q_{1}}(f) - 2\nu_{1}\right)^{2}\right)
\\
&& + \mathbf{E}\left(\left(q_{2}^{-1}\left\|DI_{q_{2}}(h)\right\|_{L^{2}(T)}^{2} - 2I_{q_{2}}(h) - 2\nu_{2}\right)^{2}\right)
\\
&& + 2\mathbf{E}\left(\left(q_{1}^{-1}\left\|DI_{q_{1}}(f)\right\|_{L^{2}(T)}^{2} - 2I_{q_{1}}(f) - 2\nu_{1}\right)\left(q_{2}^{-1}\left\|DI_{q_{2}}(h)\right\|_{L^{2}(T)}^{2} - 2I_{q_{2}}(h) - 2\nu_{2}\right)\right).
\end{eqnarray*}
We use Lemma \ref{lemmeIndep} to write
\begin{eqnarray*}
 \mathbf{E}\left(\left(q_{1}^{-1}\left\|DI_{q_{1}}(f)\right\|_{L^{2}(T)}^{2} - 2I_{q_{1}}(f) - 2\nu_{1}\right)\left(q_{2}^{-1}\left\|DI_{q_{2}}(h)\right\|_{L^{2}(T)}^{2} - 2I_{q_{2}}(h) - 2\nu_{2}\right)\right)
&=& 0.
\end{eqnarray*}
Thus,
\begin{eqnarray*}
 \mathbf{E}\left(\left(2\nu + 2Z - G_{Z}\right)^{2}\right)
&& = q_{1}^{-1}\mathbf{E}\left(\left(\left\|DI_{q_{1}}(f)\right\|_{L^{2}(T)}^{2} - 2q_{1}I_{q_{1}}(f) - 2q_{1}\nu_{1}\right)^{2}\right)
\\
&& + q_{2}^{-1}\mathbf{E}\left(\left(\left\|DI_{q_{2}}(h)\right\|_{L^{2}(T)}^{2} - 2q_{2}I_{q_{2}}(h) - 2q_{2}\nu_{2}\right)^{2}\right).
\end{eqnarray*}
Relation (\ref{equiv4}) and the calculations contained in \cite{NoPe2} imply that the above two summands vanish.
\\
\noindent It finally follows from this that
\begin{eqnarray*}
\mathbf{E}\left(\left(2\nu + 2Z - G_{Z}\right)^{2}\right) =0.
\end{eqnarray*}
\\
\textit{ Proof of $(i) \rightarrow (ii)$}. Suppose that $(ii)$ holds.  We have proven that $$\mathbf{E}\left(\left(2\nu + 2Z - G_{Z}\right)^{2}\right) =0 \Rightarrow \left\{
\begin{array}{ll}
\mathbf{E}\left(\left(\left\|DI_{q_{1}}(f)\right\|_{L^{2}(T)}^{2} - 2q_{1}I_{q_{1}}(f) - 2q_{1}\nu_{1}\right)^{2}\right) = 0 \\
\mathbf{E}\left(\left(\left\|DI_{q_{2}}(h)\right\|_{L^{2}(T)}^{2} - 2q_{2}I_{q_{2}}(h) - 2q_{2}\nu_{2}\right)^{2}\right) = 0.
\end{array}
\right.
$$ From Theorem 1.2 in \cite{NoPe2} it follows that $I_{q_{1}}(f) \sim F(\nu_{1})$ and $I_{q_{2}}(h) \sim F(\nu_{2})$. $I_{q_{1}}(f)$ and $I_{q_{2}}(h)$ being independent, we use the convolution property of Gamma random variables to state that $Z = I_{q_{1}}(f) + I_{q_{2}}(h) \sim F(\nu_{1}+\nu_{2}) \sim F(\nu)$. \qed

\begin{remark}
The proof of the above theorem shows that the affirmations (i) and (ii) are equivalent with relations (\ref{equiv2}), (\ref{equiv1}), (\ref{equiv3}) and (\ref{equiv4}).
\end{remark}
\noindent Following exactly the lines  of the proof of Theorem \ref{theoGeneral} it is possible to characterize random variables given by a sum of independent multiple integrals  that converges in law to a Gamma distribution.

\begin{theorem}
\label{Asym}
\noindent Fix $\nu_{1},\nu_{2},\nu > 0$ such that $\nu_{1}+\nu_{2}=\nu$ and let $F(\nu)$ be a real-valued random variable with characteristic function given by (\ref{fc}).
\noindent Fix two even integers $q_{1} \geq 2$ and $q_{2} \geq 2$. For any sequence  $(f_{k})_{k\geq 1} \subset L^{2}(T^{ q_{1}})$ and $(h_{k})_{k\geq 1} \subset L^{2}(T^{ q_{2}})$ ($f_{k}$ and $h_{k}$ are symmetric for every $k\geq 1$) such that
\begin{eqnarray*}
\mathbf{E}\left(I_{q_{1}}(f_{k})^{2}\right) = q_{1}!\left\|f_{k}\right\|^2_{L^{2}(T^{q_{1}})} \underset{k \rightarrow +\infty}{\longrightarrow} 2\nu_{1} \mbox{\  \  \  and\  \  \  } \mathbf{E}\left(I_{q_{2}}(h_{k})^{2}\right) = q_{2}!\left\|h_{k}\right\|^2_{L^{2}(T^{ q_{2}})}\underset{k \rightarrow +\infty}{\longrightarrow} 2\nu_{2},
\end{eqnarray*}
and such that $X_{k} = I_{q_{1}}(f_{k})$ and $Y_{k} = I_{q_{2}}(h_{k})$ are independent for any $k\geq 1$, define the random variable
\begin{eqnarray*}
Z _{k}= X_{k} + Y_{k} = I_{q_{1}}(f_{k}) + I_{q_{2}}(h_{k}) \hskip0.5cm \forall k\geq 1.
\end{eqnarray*}
Under those conditions, the following two conditions are equivalent:
\begin{enumerate}[(i)]
\item $\mathbf{E}\left(\left(2\nu + 2Z _{k}- \left\langle DZ_{k},-DL^{-1}Z_{k}\right\rangle_{L^{2}(T)}\right)^{2}\right) \underset{k \rightarrow +\infty}{\longrightarrow} 0$;
\item $Z_{k} \stackrel{\rm Law}{\underset{k \rightarrow +\infty}{\longrightarrow}}F(\nu)$;
\end{enumerate}
\end{theorem}
\noindent Cram\'er's theorem for Gamma random variables in the setting of multiple stochastic integrals is a corollary of Theorem \ref{theoGeneral}. We have the following :
\begin{theorem}
\label{cramer}
\noindent Let $Z = X + Y = I_{q_{1}}(f) + I_{q_{2}}(h)$,  $q_{1}, q_{2}\geq 2$, $f\in L^{2}(T^{q_{1}}), h\in L^{2}(T^{q_{2}})$ symmetric,  be such that $X , Y$ are independent  and
$$
\mathbf{E}\left(Z^{2}\right) = 2\nu, 
\mathbf{E}\left(X^{2}\right) = q_{1}!\left\|f\right\|^2_{L^{2}(T^{q_{1}})}=2\nu_{1}, 
\mathbf{E}\left(Y^{2}\right) = q_{2}!\left\|h\right\|^2_{L^{2}(T^{q_{2}})}=2\nu_{2}
$$
with $\nu = \nu_{1} + \nu_{2}$. Furthermore, let's assume that $Z \sim F(\nu)$. Then,
\begin{eqnarray*}
X \sim F(\nu_{1}) \mbox{\  \  \  and\  \  \  } Y \sim F(\nu_{2}).
\end{eqnarray*}
\end{theorem}
{\bf Proof: }
Theorem \ref{theoGeneral} states that $Z \sim F(\nu) \Leftrightarrow \mathbf{E}\left(\left(2\nu + 2Z - G_{Z}\right)^{2}\right) =0$ and we proved that
\begin{eqnarray*}
\mathbf{E}\left(\left(2\nu + 2Z - G_{Z}\right)^{2}\right) = \mathbf{E}\left(\left(2\nu_{1} + 2X - G_{X}\right)^{2}\right) + \mathbf{E}\left(\left(2\nu_{2} + 2Y - G_{Y}\right)^{2}\right).
\end{eqnarray*}
Both summands being positive, it follows that $\mathbf{E}\left(\left(2\nu_{1} + 2X - G_{X}\right)^{2}\right) = 0$ and $\mathbf{E}\left(\left(2\nu_{2} + 2Y - G_{Y}\right)^{2}\right) = 0$. Applying theorem \ref{theoGeneral} to $X$ and $Y$ separately gives us
$\mathbf{E}\left(\left(2\nu_{1} + 2X - G_{X}\right)^{2}\right) \Leftrightarrow X \sim F(\nu_{1})$
and $
\mathbf{E}\left(\left(2\nu_{2} + 2Y - G_{Y}\right)^{2}\right) \Leftrightarrow Y \sim F(\nu_{2}).$\qed

\vskip0.2cm

\noindent It is immediate to give an asymptotic version of Theorem \ref{cramer}.

\begin{theorem}
\label{cramerAsym}
\noindent Let $Z_{k} = X_{k} + Y_{k} = I_{q_{1}}(f_{k}) + I_{q_{2}}(h_{k})$, $f_{k}\in L^{2}(T^{q_{1}}), h_{k}\in L^{2}(T^{q_{2}})$ symmetric for  $k\geq 1$,  $q_{1}, q_{2}\geq 2$,  be such that $X_{k}, Y_{k}$ are independent for every $k\geq 1$ and
$$
\mathbf{E}\left(Z_{k}^{2}\right) \underset{k \rightarrow +\infty}{\longrightarrow} 2\nu, 
\mathbf{E}\left(X_{k}^{2}\right) = q_{1}!\left\|f\right\|^2_{L^{2}(T^{ q_{1}})}\underset{k \rightarrow +\infty}{\longrightarrow}2\nu_{1}, 
\mathbf{E}\left(Y_{k}^{2}\right) = q_{2}!\left\|h\right\|^2_{L^{2}(T^{q_{2}})} \underset{k \rightarrow +\infty}{\longrightarrow}2\nu_{2}
$$
with $\nu = \nu_{1} + \nu_{2}$. Furthermore, let's assume that $Z_{k} \underset{k \rightarrow +\infty}{\longrightarrow} F(\nu)$  in distribution. Then,
\begin{eqnarray*}
X_{k} \underset{k \rightarrow +\infty}{\longrightarrow}  F(\nu_{1}) \mbox{\  \  \  and\  \  \  } Y _{k}\underset{k \rightarrow +\infty}{\longrightarrow} F(\nu_{2}).
\end{eqnarray*}
\end{theorem}

\begin{remark} i) From Corollary  4.4. in \cite{NoPe2} it follows that actually there are no Gamma distributed random variables in a chaos of order bigger or equal than 4. (We actually conjecture that a Gamma distributed random variable given by a multiple integral can only live in the second Wiener chaos). In this sense Theorem \ref{cramer} contains a limited number of examples. By contrary, the asymptotic Cram\'er theorem (Theorem \ref{cramerAsym}) is more interesting and more general since there exists a large class of variables which are asymptotically  Gamma distributed.
\\
\noindent ii) Theorem \ref{cramer} cannot be applied directly to random variables with law $\Gamma (a, \lambda)$ (as defined in the Introduction) because  such random variables are not centered and then they cannot live in a finite Wiener chaos. But, it is not difficult to understand that if $X=I_{q_{1}}+c $ is a random variable  which is independent of $Y= I_{q_{2}}+d$ (and assume that the first two moments of $X$ and $Y$ are the same as the moment of the corresponding Gamma distributions),   and if $X+Y\sim \Gamma (a+b,\lambda) $  then $X$ has the distribution $\Gamma (a, \lambda)$ and $Y$ has the distribution $\Gamma (b, \lambda)$.
\\
\noindent iii) Several results of the paper (Lemmas 1 and 2) holds for strongly independent random variables.  Nevertheless, the key results (Theorems \ref{theoGeneral} and \ref{Asym} that allows to prove Cram\'er's theorem and its asymptotic variant are not true for strongly independent random variables (actually the implication $ii)\to i)$ in these results, whose proof is based on the differential equation satisfied by the characteristic function of the Gamma distribution, does not work.
\end{remark}

\section{Counterexample in the general case}

We will see in this section that Theorem \ref{cramer} does not hold for random variables which have a chaos decomposition into an infinite sum of multiple stochastic integrals. We construct a counterexample in this sense. What is more interesting is that the random variables defined in the below example are not only independent, they are {\it strongly independent} (see the definition above).

\begin{example}
\label{counter}
Let $\epsilon(\lambda)$ denote the exponential distribution with parameter $\lambda$ and let $b(p)$ denote the Bernoulli distribution with parameter $p$. Let $X = A -1$ and $Y = 2\varpi B -1$, where $A \sim \epsilon(1)$, $B \sim \epsilon(1)$, $\varpi \sim b(\frac{1}{2})$ and $A$, $B$ and $\varpi$ are mutually independent. This implies that $X$ and $Y$ are independent. We have $\mathbf{E}(X) = \mathbf{E}(Y) = 0$ as well as  $\mathbf{E}(X^{2}) = 1$ and $\mathbf{E}(Y^{2}) = 3$. Consider also $Z = X+Y$. Observe that $X$,$Y$ and $Z$ match every condition of theorem \ref{cramer}, but $X$ and $Y$ are not multiple stochastic integrals in a fixed Wiener chaos (see the next proposition for more details). We have the following : $Z \sim F(2)$, but $Y$ is not Gamma distributed.
\end{example}
{\bf Proof: }
We know that
\begin{eqnarray*}
\mathbf{E}\left(e^{it X}\right) = \mathbf{E}\left(e^{it (A-1)}\right) = e^{-it}\mathbf{E}\left(e^{it A}\right) = \frac{e^{-it}}{1-it}
\end{eqnarray*}
and that
\begin{eqnarray*}
\mathbf{E}\left(e^{it Y}\right) &=& \mathbf{E}\left(e^{it (2\varpi B-1)}\right) = e^{-it}\mathbf{E}\left(e^{it 2\varpi B}\right) = e^{-it}\left(\frac{1}{2}\mathbf{E}\left(e^{it 2B}\right) + \frac{1}{2}\right)
\\
&=& e^{-it}\left(\frac{1}{2}\frac{1}{1-2it} + \frac{1}{2}\right) = e^{-it}\frac{1-it}{1-2it}.
\end{eqnarray*}
Observe at this point that the characteristic function of $Y$ proves that $Y$ is not Gamma distributed.
Let us compute the characteristic function of $Z$. We have
\begin{eqnarray*}
\mathbf{E}\left(e^{it Z}\right) = \mathbf{E}\left(e^{it (X+Y)}\right) = \mathbf{E}\left(e^{it X}\right)\mathbf{E}\left(e^{it Y}\right) = \frac{e^{-it}}{1-it}e^{-it}\frac{1-it}{1-2it} = \frac{e^{-2it}}{1-2it} = \mathbf{E}\left(e^{it F(2)}\right).
\end{eqnarray*}
\qed

\begin{remark}
It is also possible to construct a   similar example for the laws $\Gamma (a, \lambda), \Gamma (b, \lambda)$ instead of $F(\nu_{1}), F(\nu_{2})$.
\end{remark}
\vskip0.2cm

\noindent The following proposition shows that this counterexample accounts for independent random variables but also for strongly independent random variables.

\vskip0.2cm

\begin{prop}
$X$ and $Y$ as defined in Example  \ref{counter} are strongly independent.
\end{prop}
{\bf Proof: }
In order to prove that $X$ and $Y$ are strongly independent, we need to compute their Wiener chaos expansions in order to emphasize the fact that all the components of these Wiener Chaos expansions are mutually independent. Consider a standard Brownian motion $B$ indexed on $L^{2}(T) = L^{2}(\left(0,T\right))$. Consider $h_{1},...,h_{5} \in L^{2}(T)$ such that $\left\|h_{i}\right\|_{L^{2}(T)} = 1$ for every $1 \leq i \leq 5$ and such that $W(h_{i})$ and $W(h_{j})$ are independent for every $1 \leq i,j \leq 5, i \neq j$. First notice that the random variables $A = \frac{1}{2}\left(W(h_{1})^{2} + W(h_{2})^{2}\right)$ and $B = \frac{1}{2}\left(W(h_{4})^{2} + W(h_{5})^{2}\right)$ are independent (this is obvious) and have the exponential distribution with parameter 1.
Also, note that the random variable $\varpi = \frac{1}{2}\mbox{sign}(W(h_{3})) + \frac{1}{2}$ has the Bernoulli distribution and is independent from $A$ and $B$. As in Example \ref{counter}, set $X = A - 1$ and $Y = 2\varpi B - 1$. $X$ and $Y$ are as defined in Example \ref{counter}. Let us now compute their Wiener chaos decompositions. We have
\begin{eqnarray*}
A &=& \frac{1}{2}\left(W(h_{1})^{2} + W(h_{2})^{2}\right) = \frac{1}{2}\left(I_{1}(h_{1})^{2} + I_{1}(h_{2})^{2}\right)
= \frac{1}{2}\left(2 + I_{2}(h_{1}^{\otimes 2}) + I_{2}(h_{2}^{\otimes 2})\right),
\end{eqnarray*}
and similarly $B = \frac{1}{2}\left(2 + I_{2}(h_{4}^{\otimes 2}) + I_{2}(h_{5}^{\otimes 2})\right).$ Therefore, we have
\begin{eqnarray*}
X = I_{2}\left(\frac{h_{1}^{\otimes 2} + h_{2}^{\otimes 2}}{2}\right).
\end{eqnarray*}
From \cite{HuNu1}, Lemma 3, we know that
\begin{eqnarray*}
\mbox{sign}(W(h_{3})) = \sum_{k \geq 0}b_{2k+1}I_{2k+1}(h_{3}^{\otimes (2k+1)}),
\end{eqnarray*}
where $b_{2k+1} = \frac{2(-1)^{k}}{(2k+1)\sqrt{2\pi}k!2^{k}}.$ It follows that $
\varpi = \frac{1}{2} + \frac{1}{2}\sum_{k \geq 0}b_{2k+1}I_{2k+1}(h_{3}^{\otimes (2k+1)}),
$
and
\begin{eqnarray*}
Y &=& (1 + \sum_{k \geq 0}b_{2k+1}I_{2k+1}(h_{3}^{\otimes (2k+1)}))(1 + \frac{1}{2}I_{2}(h_{4}^{\otimes 2}) + \frac{1}{2}I_{2}(h_{5}^{\otimes 2})) - 1
\\
&=& \frac{1}{2}I_{2}(h_{4}^{\otimes 2}) + \frac{1}{2}I_{2}(h_{5}^{\otimes 2})
 + \sum_{k \geq 0}b_{2k+1}I_{2k+1}(h_{3}^{\otimes (2k+1)})
+ \frac{1}{2}\sum_{k \geq 0}b_{2k+1}I_{2k+1}(h_{3}^{\otimes (2k+1)})I_{2}(h_{4}^{\otimes 2})
\\
&& + \frac{1}{2}\sum_{k \geq 0}b_{2k+1}I_{2k+1}(h_{3}^{\otimes (2k+1)})I_{2}(h_{5}^{\otimes 2}).
\end{eqnarray*}
Using the multiplication formula for multiple stochastic integrals, we obtain
\begin{eqnarray*}
Y &=& \frac{1}{2}I_{2}(h_{4}^{\otimes 2}) + \frac{1}{2}I_{2}(h_{5}^{\otimes 2})+ \sum_{k \geq 0}b_{2k+1}I_{2k+1}(h_{3}^{\otimes (2k+1)})
\\
&& +  \frac{1}{2}\sum_{k \geq 0}b_{2k+1}\sum_{r=0}^{(2k+1) \wedge 2}r!\binom{2}{r}\binom{2k+1}{r}I_{2k+3-2r}(h_{3}^{\otimes (2k+1)} \otimes_{r} h_{4}^{\otimes 2})
\\
&& +  \frac{1}{2}\sum_{k \geq 0}b_{2k+1}\sum_{r=0}^{(2k+1) \wedge 2}r!\binom{2}{r}\binom{2k+1}{r}I_{2k+3-2r}(h_{3}^{\otimes (2k+1)} \otimes_{r} h_{5}^{\otimes 2}).
\end{eqnarray*}
At this point, it is clear that $X$ and $Y$ are strongly independent.\qed

\end{document}